\documentclass[a4paper,12pt]{article}
\usepackage{amsfonts}
\usepackage{amssymb,amsmath} 
\numberwithin{equation}{section}
\newtheorem{theorem}{Theorem}[section]

\newtheorem{proposition}{Proposition}[section]

\usepackage{cite}

\begin{document}
\thispagestyle{empty}

\vbox{\hbox{\footnotesize Bulletin of the \textit{Transilvania} University of Bra\c sov $\bullet$ Vol 11(60), No. 1 - 2018}}
\vbox{\hbox{\footnotesize Series III: Mathematics, Informatics, Physics, xx-xx}}
\vskip 1.5truecm

\begin{center}
{\large{\bf \textbf{$\mathbb{R}$-COMPLEX FINSLER SPACES WITH INFINITE SERIES $(\alpha, \beta)$-METRIC}}}
\end{center}
\medskip

\begin{center}
{\bf GAUREE SHANKER \footnote{Department of Mathematics and Statistics, \textit{Central} University of Punjab, Bhatinda-151 001, Punjab, India, e-mail: grshnkr2007@gmail.com} and RUCHI KAUSHIK SHARMA\footnote{Department of Mathematics and Statistics, \textit{Banasthali} University, Banasthali-304022, Rajasthan, India, e-mail: ruchikaushik07@gmail.com}}
\end{center}

\pagestyle{myheadings}
\markboth{GAUREE SHANKER and RUCHI KAUSHIK SHARMA
}{$\mathbb{R}$-Complex Finsler Space with Infinite Series Metric}

\bigskip

\begin{abstract}
In the present paper, the notion of $\mathbb{R}$-complex Finsler space with Infinite Series ($\alpha, \beta$)- metric $\dfrac{\beta^2}{\beta - \alpha}$ is defined. The Fundamental metric fields $g_{ij}$, $g_{i\bar{j}}$, their determinants and the inverse of these tensor fields are obtained. Also some properties of these spaces are studied.
\end{abstract}

\bigskip
\textbf{2010 AMS Classification:} 53B40, 53C60\\
\textbf{Keywords}: Complex Finsler space, $\mathbb{R}$-complex Finsler space, Infinite Series metric

\section{Introduction}
A Finsler metric L($\alpha, \beta$) in a smooth manifold $M^n$ is called an ($\alpha, \beta$)- metric, if L is a positively homogeneous function of degree one of a Riemannian metric $\alpha$ = $\sqrt{a_{ij}y^iy^j}$ and a one-form $\beta$ = $b_i(x)y^i$ on $M^n$. The interesting and important examples of an ($\alpha, \beta$)-metric are Randers metric $\alpha + \beta$, Kropina metric $\dfrac{\alpha^2}{\beta}$ and Matsumoto-metric $\dfrac{\alpha^2}{\alpha - \beta}$. \\

In the r-th series $(\alpha, \beta)$-metric $L(\alpha, \beta)$ = $\beta \sum\limits_{k=0}^{r}$ $(\dfrac{\alpha}{\beta})^k$, where we assume $\alpha < \beta$. If r = $\infty$, then this metric is expressed in the form $L(\alpha, \beta)$ = $\dfrac{\beta^2}{\beta - \alpha}$, and is called an infinite series $(\alpha, \beta)$-metric. Interestingly this metric is the difference between a Randers metric and a Matsumoto metric. Recently, some authors 
([13], [14], [15] and [16]) have worked on Finsler spaces with infinite series ($\alpha, \beta$)-metric and have produced some important and interesting results. \\

The notion of complex Finsler space appeared for the first time in a paper written by Rizza in 1963, [12] , as a generalization of the similar notion from the real case, requiring the homogeneity of the fundamental function with respect to the fibre variables, for any complex scalars $\lambda$. The first example comes from the complex hyperbolic geometry and was given by S. Kobayashi in 1975, [9]. The Kobayashi metric has given an impulse to the study of complex Finsler geometry. \\

A complex Finsler geometry, which contains many interesting results, has been developed in the papers ([1], [3], [4] and [7] etc). In paper [10], the well-known definition of a complex Finsler space ([1], [3]) was extended, reducing the scalars to $\lambda \in \mathbb{R}$. The outcome was a new class of Finsler space called the $\mathbb{R}$-complex Finsler spaces. Some important results on $\mathbb{R}$-Complex Finsler spaces have been obtained in ( [2], [7]). \\

Our interest in this class of Finsler spaces issues from the fact that the Finsler geometry means, first of all, distance and this refers to curves depending on the real parameter. \\

In the present paper following the ideas from real Finsler spaces with Infinite series metrics, we introduce the similar notions on $\mathbb{R}$-complex Finsler spaces with infinite series ($\alpha, \beta$)-metric.

\section{$\mathbb{R}$-Complex Finsler spaces}
Let M be a complex manifold with $dim_c$M = n, ($z^k$) be local complex coordinates in a chart (U, $\phi$) and $T'M$ its holomorphic tangent bundle. It has a natural structure of complex manifold, $dim_cT'M$ = 2n and the induced coordinates in a local chart on u $\in T'M$ are denoted by u = ($z^k$, $\eta^k$). The changes of local coordinates in u are given by the rules:
\begin{equation}\label{eq1}
z'^{k} = z'^{k}(z);  \eta'^{k} = \frac{\partial z'^{k}}{\partial z^j}\eta^j. 
\end{equation}

The natural frame $\displaystyle  \left\lbrace \frac{\partial}{\partial z^k}, \frac{\partial}{\partial \eta^k}\right\rbrace $ of $T'_u(T'M)$ transforms with the Jacobi matrix of equation \eqref{eq1} changes,\\

\begin{center}
 $\displaystyle \dfrac{\partial}{\partial z^k} = \dfrac{\partial z'^{j}}{\partial z^k}\dfrac{\partial}{\partial z'^{j}} + \dfrac{\partial^2 z'^{j}}{\partial z^{k}\partial z^h}\eta^h \dfrac{\partial}{\partial \eta'^{j}};5
  \dfrac{\partial}{\partial \eta^k} = \dfrac{\partial z'^{j}}{\partial z^k}\dfrac{\partial}{\partial \eta'^{j}}$.\\
\end{center}

A complex non-linear connection, briefly (c.n.c), is a supplementary distribution H($T'M$) to the vertical distribution V($T'M$) in $T'(T'M)$. The vertical distribution is spanned by $\displaystyle \left\{\frac{\partial}{\partial \eta^k} \right\}$ and an adapted frame in H($T'M$) is $ \displaystyle \frac{\delta}{\delta z^k} = \frac{\partial}{\partial z^k} - N^j_k \frac{\partial}{\partial \eta^j}$, where $N^j_k$ are the coefficients of the (c.n.c.) and they have a certain rule of change at (2.1) so that $\displaystyle \frac{\delta}{\delta z^k}$ transform like vectors on the base manifold M. Next we use the abbreviations: $\displaystyle \partial_k = \frac{\partial}{\partial z^k}, \delta_k = \frac{\delta}{\delta z^k}, \dot{\partial}_k = \frac{\partial}{\partial \eta^k}$ and $\partial_{\bar{k}},\dot{\partial}_{\bar{k}}, \delta_{\bar{k}}$ for their conjugates. The dual adapted basis of $\displaystyle \{\delta_k, \dot{\partial}_k\}$ are $\displaystyle \{dz^k, \delta \eta^k = d\eta^k + N^k_j dz^j\}$ and $\{d\bar{z}^k, \delta \bar{\eta}^k\}$ their conjugates. \\

We recall that the homogeneity of the metric function of a complex Finsler space is with respect to all complex scalars and the metric tensor of the space, is a Hermitian one [11]. \\

An $\mathbb{R}$-complex Finsler metric on M is a continuous function F : $T'M \rightarrow \mathbb{R}_+ $ satisfying:
\begin{enumerate}
\item[(i)] L := $F^2$ is smooth on ${T'M}$ (except the 0 sections);
\item[(ii)] $F(z, \eta) \geq 0$, the equality holds if and only if $\eta$ = 0;
\item[(iii)] $F( z, \lambda\eta, \bar{z}, \lambda \bar{\eta}) = \lvert \lambda \rvert F(z, \eta, \bar{z}, \bar{\eta}), \forall \lambda \in \mathbb{R}$.\\
\end{enumerate}

It follows that L is (2, 0) homogeneous with respect to the real scalars $\lambda$, and Purcaru [11] proved that the following identities are fulfilled: 
\begin{equation}\label{eq2.2}
\dfrac{\partial L}{\partial \eta^i}\eta^i + \dfrac{\partial L}{\partial \bar{\eta}^i}\bar{\eta}^i = 2L;   g_{ij}\eta^i + g_{\bar{j}i}\bar{\eta}^i = \dfrac{\partial L}{\partial \eta^j};
\end{equation}

\begin{equation}\label{eq2.3}
 \dfrac{\partial g_{ik}}{\partial \eta^i}\eta^j + \dfrac{\partial g_{ij}}{\partial \bar{\eta}^j}\bar{\eta}^j = 0; 
\dfrac{\partial g_{i\bar{k}}}{\partial \eta^j}\eta^j + \dfrac{\partial g_{i\bar{k}}}{\partial \bar{\eta}^j}\bar{\eta}^j = 0,
\end{equation}

\begin{equation}\label{eq2.4}
2L = g_{ij}\eta^i \eta^j + g_{\bar{i}\bar{j}} \bar{\eta}^i \bar{\eta}^j + 2g_{i\bar{j}} \eta^i \bar{\eta}^j,
\end{equation}
where
\begin{equation*}
g_{ij} = \dfrac{\partial^{2}L}{\eta^i\eta^j}, g_{i\bar{j}} = \dfrac{\partial^{2}L}{\eta^i\bar{\eta}^j};  g_{\bar{i}\bar{j}} := \dfrac{\partial^{2}L}{\bar{\eta}^i\bar{\eta}^j}
\end{equation*}
are the metric tensors of space.
\section{$\mathbb{R}$ - complex Finsler space with infinite series ($\alpha, \beta$)-metric}
 
An $\mathbb{R}$ - complex Finsler space (M, F) is called $\mathbb{R}$ - complex Infinite series space if
 \begin{equation}\label{eq3.1}
      F(\alpha, \beta) = \dfrac{\beta^2}{\beta-\alpha}, \alpha \neq \beta, 
      \end{equation}
       where 
      \begin{equation*}
      \alpha^2(z, \eta, \bar{z}, \bar{\eta}) = Re\{a_{ij}\eta^i\eta^j + a_{i\bar{j}} \eta^i\bar{\eta}^{j}\},
      \end{equation*}
      \begin{equation*} 
       \beta(z, \eta, \bar{z}, \bar{\eta}) = Re\{b_i\eta^i\}
      \end{equation*}
      with $ a_{ij} = a_{ij}(z), a_{i\bar{j}}$ and $b = b_i(z)dz^i$ is a (1, 0) - differential form. \\
 
 \noindent \textbf{Remark.} L = $F^2$ is $\mathbb{R}$-complex Finsler space with $(\alpha, \beta)$-metric. \\
      
\noindent  Taking into account the 2-homogeneity condition of L: 
            \begin{equation*}
            L(\alpha(z, \lambda\eta, \bar{z}, \lambda\bar{\eta}), \beta(z, \lambda\eta, \bar{z}, \lambda\bar{\eta})) = \lambda^2(\alpha(z, \eta, \bar{z}, \bar{\eta}), \beta(z, \eta, \bar{z}, \bar{\eta})), \lambda\in \mathbb{R}_+. 
      \end{equation*}
      
      \begin{proposition}
      In a $\mathbb{R}$-complex Finsler space with ($\alpha, \beta$)-metric the following equalities hold:
      \begin{equation*}
     \alpha L_{\alpha} + \beta L_{\beta} = 2L,
     \end{equation*}
         \begin{equation*}
                \alpha L_{\alpha \alpha} + \beta L_{\alpha \beta} = L_{\alpha},
                                                          \end{equation*}
                                                          \begin{equation*} 
                                                          \alpha L_{\alpha \beta} + \beta L_{\beta \beta} = L_{\beta}, 
                                                          \end{equation*}
                                                          \begin{equation*}
                                                          \alpha^2 L_{\alpha \alpha} + 2 \alpha \beta L_{\alpha \beta} + \beta^2 L_{\beta \beta} = 2L,
                                                          \end{equation*}
                                                          where
  \begin{equation*}
  L_{\alpha} = \dfrac{\partial L}{\partial \alpha}, L_{\beta} = \dfrac{\partial L}{\partial \beta}, L_{\alpha \alpha} = \dfrac{\partial^2L}{\partial \alpha^2}, L_{\beta \beta} = \dfrac{\partial^2 L}{\partial \beta^2}, L_{\alpha \beta} = \dfrac{\partial^2L}{\partial \alpha \partial \beta}.
  \end{equation*}
  \end{proposition}
    \textbf{Particular case.} For an $\mathbb{R}$-complex Finsler space with special $(\alpha, \beta)$-metric  L($\alpha, \beta$) = $\dfrac{\beta^2}{\beta-\alpha}$,                                  
      we have:
 \begin{equation}
      L_\alpha = \dfrac{2\beta^4}{(\beta - \alpha)^3},
                \end{equation}
      \begin{equation}
                     L_\beta  = \dfrac{2\beta^3 (\beta - 2\alpha)}{(\beta - \alpha)^3},
                               \end{equation}
       \begin{equation}
       L_{\alpha \alpha} = \dfrac{6 \beta^4}{(\beta - \alpha)^4},
       \end{equation} 
       \begin{equation}
       L_{\alpha \beta} = \dfrac{2\beta^3 (\beta - 4\alpha)}{(\beta - \alpha)^4},
       \end{equation}
       \begin{equation}
       L_{\beta \beta} = \dfrac{2\beta^2(\beta^2 - 4\alpha \beta + 6\alpha^2)}{(\beta - \alpha)^4}.
       \end{equation} 
        \begin{eqnarray}
                                                   \alpha L_\alpha + \beta L_\beta &=& \alpha \Bigl[\dfrac{2\beta^4}{(\beta - \alpha)^3}\Bigr] + \beta  \Bigl[2F{\dfrac{(2\beta^4 - 4\alpha \beta^3)}{(\beta - \alpha)^3}}\Bigr]
                                                   \nonumber \\
                                                   &&
                                                   = \dfrac{2\beta^5 - 2\alpha\beta^4}{(\beta - \alpha)^3} = 2\dfrac{\beta^4}{(\beta - \alpha)^2} = 2L.
                                                   \end{eqnarray}
        \begin{eqnarray}
       \alpha L_{\alpha \alpha} + \beta L_{\alpha \beta} &=& \alpha \Bigl[\dfrac{6\beta^4}{(\beta - \alpha)^4}\Bigr] + \beta \Bigl[\dfrac{2\beta^4 - 8\alpha\beta^3}{(\beta - \alpha)^4}\Bigr] \nonumber \\
       &&
       = \dfrac{6\alpha\beta^4 + 2\beta^5 - 8\alpha\beta^4}{(\beta - \alpha)^4} = \dfrac{2\beta^4}{(\beta - \alpha)^3} = L_{\alpha}.
       \end{eqnarray}
       Similarly, we can prove other equalities.
       \par
       
       In the following we propose to determine the metric tensors of an $\mathbb{R}$-complex Finsler space with Infinite series metric, i.e., 
                \begin{equation}
                 g_{ij} = \dfrac{\partial^{2}L(z, \eta, \bar{z}, \lambda\bar{\eta})}{\eta^i\eta^j}, g_{i\bar{j}} = \dfrac{\partial^{2}L(z, \eta, \bar{z}, \lambda\bar{\eta})}{\eta^i\bar{\eta}^j}.
               \end{equation}
               Each of these being of interest in the following:
              
               \par
               We consider:
               \begin{equation*}
               \dfrac{\partial \alpha}{\partial \eta^i} = \dfrac{1}{2\alpha}(a_{ij}\eta^j + a_{i\bar{j}}\bar{\eta}^j) = \dfrac{1}{2\alpha}l_i, \dfrac{\partial \beta}{\partial \eta^i} = \dfrac{1}{2}b_i,
               \end{equation*}
               \begin{equation*}
                       \dfrac{\partial \alpha}{\partial \bar{\eta}^i} = \dfrac{1}{2\alpha}(a_{\bar{i}\bar{j}}\bar{\eta}^j + a_{i\bar{j}}{\eta}^j) = 
                       \dfrac{1}{2\alpha} l_{\bar{i}}, \dfrac{\partial \beta}{\partial \bar{\eta}^i} = \dfrac{1}{2}b_{\bar{i}},
              \end{equation*}
              where, 
              \begin{equation*}
              l_i = a_{\bar{i}\bar{j}}\bar{\eta}^j + a_{i\bar{j}}{\eta}^j, l_{\bar{j}} = a_{\bar{i}\bar{j}}\bar{\eta}^i + a_{i\bar{j}}{\eta}^i.
              \end{equation*}
              We find immediately: 
              \begin{equation*}
              l_i \eta^i + l_{\bar{j}} \bar{\eta}^j = 2 \alpha^2.
              \end{equation*}
              We denote:
              \begin{equation*}
              \eta_i = \dfrac{\partial L}{\partial \eta^i} = \dfrac{\partial }{\partial \eta^i}{F^2} = 2F \dfrac{\partial }{\partial \eta^i} \Bigl(\dfrac{\beta^2}{\beta - \alpha}\Bigr),
              \end{equation*}
              \begin{equation*}
              \eta_i = \rho_0 l_i + \rho_1 b_i,
              \end{equation*}
              where
              \begin{equation*}
              \rho_0 = \dfrac{1}{2}\alpha^{-1}L_{\alpha} = \dfrac{\beta^4}{\alpha(\beta - \alpha)^3},
              \end{equation*}
              and
              \begin{equation*}
              \rho_1 = \dfrac{1}{2}L_{\beta} = \dfrac{\beta^3(\beta - 2\alpha)}{(\beta - \alpha)^3}.
              \end{equation*}
              Differentiating $\rho_0$ and $\rho_1$ w.r.t. $\eta ^j$ and $\bar{\eta}^j$ respectively, we obtain:
              \begin{align*}
              \dfrac{\partial \rho_0}{\partial \eta^j} = \dfrac{\partial}{\partial \eta^j} \Bigl[\dfrac{\beta^4}{\alpha(\beta - \alpha)^3}\Bigr]
              = \dfrac{\beta^4 (4\alpha - \beta)}{2\alpha^3(\beta - \alpha)^4}l_j + \dfrac{\beta^3 (\beta - 4\alpha)}{2\alpha(\beta - \alpha)^4}b_j
               = \rho_{-2} l_j + \rho_{-1} b_j.
              \end{align*}
              and
              \begin{align*}
                     \dfrac{\partial \rho_0}{\partial \bar{\eta}^j} = \dfrac{\partial}{\partial \bar{\eta}^j} \Bigl[\dfrac{(\alpha^3 - 2 \alpha^2 \beta)}{(\alpha - \beta)^3}\Bigr]
                     = \dfrac{\beta^4 (4\alpha - \beta)}{2\alpha^3(\beta - \alpha)^4}l_{\bar{j}} + \dfrac{\beta^3 (\beta - 4\alpha)}{2\alpha(\beta - \alpha)^4}b_{\bar{j}}
                     = \rho_{-2} l_{\bar{j}} + \rho_{-1} b_{\bar{j}}.
                     \end{align*}
                     Similarly
                 \begin{equation*}
                 \dfrac{\partial \rho_1}{\partial \eta^i} = \rho_{-1} l_i + \mu_0 b_i, \dfrac{\partial \rho_1}{\partial \bar{\eta}^i} = \rho_{-1}l_{\bar{i}} + \mu_0 b_{\bar{i}},
                 \end{equation*}
                 where
                 \begin{equation*}
                 \rho_{-2} = \dfrac{\alpha L_{\alpha \alpha} - L_\alpha}{4 \alpha^3}, \rho_{-1} = \dfrac{L_{\alpha \beta}}{4\alpha}, \mu_0 = \dfrac{L_{\beta \beta}}{4}.
                           \end{equation*}
                           \begin{proposition} 
                           The invariants of an $\mathbb{R}$ - complex Finsler space with Infinite series metric:\\
                           $\rho_0, \rho_1$, $\rho_{-2}$, $\rho_{-1}$ and $\mu_0$ are given by:
                           \begin{equation}
                           \begin{split}
                           \rho_0 &= \dfrac{1}{2}\alpha^{-1} L_{\alpha} = \dfrac{\beta^4}{\alpha(\beta - \alpha)^3},\\ 
                           \rho_1 &= \dfrac{1}{2}L_{\beta} = \dfrac{\beta^3(\beta - 2\alpha)}{(\beta - \alpha)^3},\\
                           \rho_{-2} &= \dfrac{\beta^4(4\alpha - \beta)}{2\alpha^3(\beta - \alpha)^4},\\
                           \rho_{-1} &= \dfrac{\beta^3(\beta - 4\alpha)}{2\alpha(\beta - \alpha)^4},\\
                           \mu_0 &= \dfrac{\beta^2(\beta^2 - 4\alpha\beta + 6\alpha^2)}{2(\beta - \alpha)^4}\nonumber.\\
                           \end{split}
                            \end{equation}
                            \end{proposition}
               \begin{theorem}
               The fundamental metric tensors of $\mathbb{R}$ - complex Finsler space with Infinite series metric are given by
               \begin{equation}
               \begin{split}
               g_{ij} &= \rho_0 a_{ij} + \rho_{-2} l_i l_j + \mu_0 b_i b_j + \rho_{-1}(b_j l_i + b_i l_j)\\
                      &= \dfrac{\beta^4}{\alpha(\beta - \alpha)^3}a_{ij} + \dfrac{\beta^4(4\alpha - \beta)}{2\alpha^3(\beta - \alpha)^4}l_i l_j + \dfrac{\beta^2(\beta^2 - 4\alpha\beta + 6\alpha^2)}{2(\beta - \alpha)^4}b_i b_j + \dfrac{\beta^3(\beta - 4\alpha)}{2\alpha(\beta - \alpha)^4}(b_j l_i + b_i l_j).
               \end{split}
               \end{equation}
               \begin{equation}
               \begin{split}
                g_{i\bar{j}} &= \dfrac{\beta^4}{\alpha(\beta - \alpha)^3}a_{i\bar{j}} + \dfrac{\beta (4\beta - \alpha)}{2(\alpha - \beta)^4}l_i l_{\bar{j}} 
                 + \dfrac{\beta^2(\beta^2 - 4\alpha\beta + 6\alpha^2)}{2(\beta - \alpha)^4}b_i b_{\bar{j}} + \dfrac{\beta^3(\beta - 4\alpha)}{2\alpha(\beta - \alpha)^4}(b_{\bar{j}} l_i + b_i l_{\bar{j}}).
                \end{split}
                 \end{equation}
                                \end{theorem}
               or, in the equivalent form:
               \begin{equation}
               g_{ij} = \rho_0[a_{ij} - \sigma_1l_il_j + \sigma_2b_ib_j + \sigma_3\eta_i\eta_j],
               \end{equation}
               \begin{equation}
               g_{i\bar{j}} = \rho_0[a_{i\bar{j}} - \sigma_1l_il_{\bar{j}} + \sigma_2b_ib_{\bar{j}} + \sigma_3\eta_i\eta_{\bar{j}}],
               \end{equation}
                where\\
                \begin{equation}
                \sigma_1 = \dfrac{(\beta - \alpha)^5(\beta - 4\alpha)}{2\alpha^2(\beta - 2\alpha)},
                \end{equation}
                \begin{equation}
                \sigma_2 = -\dfrac{\alpha^3}{\beta^2(\beta - \alpha)},
                \end{equation}
                \begin{equation}
                \sigma_3 = \dfrac{\alpha(\beta - \alpha)^5(\beta - 4\alpha)}{2\beta^8(\beta - 2\alpha)}.
                \end{equation}
             
                \textit{Proof.} Using the relation (30) in Theorem 2.1 given in [2] by direct calculations we obtain the results.
                \par 
                The next objective is to obtain the inverse and determinant of the tensor field $g_{ij}$.
                 \begin{proposition}
                            Suppose:\\
                            $\bullet$ $(Q_{ij})$is a non-singular n $\times$ n complex matrix with inverse $(Q^{ji})$;\\
                            $\bullet$ $C_i$ and $C_{\bar{i}} := \bar{C_i}$, i = 1,..., n are complex numbers;\\
                            $\bullet$ $C^i := Q^{ji}C_j$ and its conjugates; $C^2 := C^iC_i = \bar{C}^iC_{\bar{i}};H_{ij} := Q_{ij} \pm C_iC_j$. \\
                            Then\\
                            (i) $det(H_{ij}) = (1 \pm C^2)det(Q_{ij})$,\\
                            (ii) Whenever $1 \pm C^2 \neq 0$, the matrix $(H_{ij})$ is invertible and in this case its inverse is $H^{ji} = Q^{ji}\mp \dfrac{1}{1\pm C^2}C^iC^j$.
                            \end{proposition}      
                 \begin{theorem}
                 For a non-Hermitian $\mathbb{R}$-complex Finsler space with special ($\alpha, \beta$)-metric L$(\alpha, \beta)$ = $\dfrac{\beta^2}{\beta - \alpha}$, $\alpha \neq \beta$ with $a_{i\bar{j}}$ = 0 we have
                  \begin{multline}
                 g^{ij} = \dfrac{1}{\rho_0}\Biggl[a^{ji} + \Bigl(\dfrac{\sigma_1}{1-\gamma \sigma_1} - \dfrac{\epsilon^2 \sigma_1^2 \sigma_2}{\tau(1 - \gamma \sigma_1)^2}\Bigr)\eta^i\eta^j - \dfrac{\sigma_2b^ib^j}{\tau} - \dfrac{\epsilon \sigma_1 \sigma_2 (b^i \eta^j + b^j \eta^i)}{\tau(1-\gamma \sigma_1)}  \\
                  - \dfrac{\Omega^2 \eta^i \eta^j + \Omega \Gamma(\eta^ib^j + \eta^jb^i) + \Gamma^2b^ib^j}{1 + (\Omega \gamma + \Gamma \epsilon)\sqrt{\sigma_3}} \Biggr],
                 \end{multline} 
                 where\\
                 \begin{eqnarray}
                 l_i &=& a_{ij}\eta^j, \gamma = a_{ij}\eta^j\eta^k = l_k \eta^k, \epsilon = b_j \eta^j, \omega = b_jb^j,\nonumber \\
                 &&
                 b^k = a^{jk}b_j, b_l = b^ka_{kl}, \delta = a_{jk}\eta^jb^k = l_kb^k, l^j = a^{ji}l_i = \eta^j.
                 \end{eqnarray}
                 \begin{multline}
                  det(g_{ij})=(\rho_0)^n [1 + (\Omega \gamma + \Gamma \epsilon)\sqrt{\sigma_3}]\Bigl[1 + \omega +\\ 
                   \frac{\sigma_1 \epsilon^2}{1 - \sigma_1 \gamma}\Bigr](1 - \sigma_1 \gamma)det(a_{ij}),
                 \end{multline}
                 where
                 \begin{eqnarray}
                 \rho_0 &=& \dfrac{1}{2}\alpha^{-1}L_{\alpha} = \dfrac{\beta^4}{\alpha(\beta - \alpha)^3},
                 \end{eqnarray}
                 \begin{eqnarray}
                 \Omega = 1 + \Bigl(\dfrac{\sigma_1}{1-\gamma \sigma_1} - \dfrac{\epsilon^2 \sigma_1^2 \sigma_2}{\tau(1 - \gamma \sigma_1)^2}\Bigr)\gamma - \Bigl[\dfrac{\sigma_1 \sigma_2}{\tau (1 - \gamma \sigma_1)}\Bigr]\epsilon^2,
                 \end{eqnarray}
                 \begin{eqnarray}
                 \Gamma = - \dfrac{\sigma_2}{\tau} \epsilon + \Bigl[\dfrac{\sigma_1 \sigma_2 \epsilon}{\tau (1 - \gamma \sigma_1)}\Bigr] \gamma,
                 \end{eqnarray}
                 \begin{equation}
                                 \sigma_1 = \dfrac{(\beta - \alpha)^5(\beta - 4\alpha)}{2\alpha^2(\beta - 2\alpha)},
                                 \end{equation}
                                 \begin{equation}
                                 \sigma_2 = -\dfrac{\alpha^3}{\beta^2(\beta - \alpha)},
                                 \end{equation}
                                 \begin{equation}
                                 \sigma_3 = \dfrac{\alpha(\beta - \alpha)^5(\beta - 4\alpha)}{2\beta^8(\beta - 2\alpha)}.
                                 \end{equation}
                 \end{theorem}
                 \textit{Proof.} To prove the claims we apply  proposition (3.3) in the following steps. We write $g_{ij}$ from 3.11 in the form:
                 \begin{equation}
                  g_{ij} = \rho_0[a_{ij} - \sigma_1l_il_j + \sigma_2b_ib_j + \sigma_3\eta_i\eta_j],
                 \end{equation}
                 \textbf{Step 1.} Set $Q_{ij}$ = $a_{ij}$ and $C_i = \sqrt{\sigma_1}l_i$. By applying Proposition (3.3) we obtain $Q^{ji} = a^{ji}$, $C^2 = C_iC^i = \sqrt{\sigma_1}l_i \times Q^{ji}\times C_j = \sqrt{\sigma_1}l_i \times a^{ji}\times \sqrt{\sigma_1}l_j = \sigma_1 \times l_ia^{ji}l_j = \sigma_1 \times a_{jk} l_i \eta^i = \sigma_1 \gamma $ and $ 1 - C^2 = (1-\sigma_1 \gamma)$.
                         \\
                         So, the matrix $H_{ij} = a_{ij} - \sigma_1 l_il_j$ is invertible with
                         \begin{equation*}
                         H^{ij} = a^{ji} + \dfrac{\sigma_1 \eta^i \eta^j}{1 - \sigma_1 \gamma},\\
                                \end{equation*}
                                \begin{equation*}
                                det(a_{ij} - \sigma_1l_il_j) = (1-\sigma \gamma)det(a_{ij}),
                                \end{equation*}
                                \textbf{Step2.} Now, take\\
                                
                                \begin{eqnarray}
                                Q_{ij} &=& a_{ij} - \sigma_1 l_il_j,\nonumber \\
                                &&
                                C_i = \sqrt{\sigma_2}b_i,\nonumber \\
                                &&
                                Q^{ji} = a^{ji} + \dfrac{\sigma_1 \eta^i\eta^j}{1 - \sigma_1 \gamma},
                                \end{eqnarray}
                                
                                and\\
                                \begin{eqnarray}
                                C^i &=& Q^{ji}C_j \nonumber \\
                                &&
                                = a^{ji} + \dfrac{\sigma_1 \eta^i \eta^j}{1 - \sigma_1 \gamma}, \nonumber \\
                                                \end{eqnarray}
                                      So,
                                     \begin{eqnarray}
                                     C^2 &=& \sigma_2 \Bigl(\omega + \dfrac{\sigma_1 \epsilon^2}{1 - \sigma_1 \gamma}\Bigr),\nonumber \\
                                     &&
                                     1 + C^2 = 1 + \sigma_2 \Bigl( \omega + \dfrac{\sigma_1 \epsilon^2}{ 1 - \sigma_1 \gamma}\Bigr) \neq 0.
                                     \end{eqnarray}
                                     Now,
                                       inverse of $ H_{ij} = a_{ij} - \sigma_1 l_il_j + \sigma_2 b_ib_j$, 
                                       exists and
                                        it is 
                                        \begin{equation}
                                        \begin{split}
                                        H^{ji} &= Q^{ji} \pm \dfrac{1}{1 + C^2}C^iC^j \\
                                         &= \Bigl(a^{ji} + \dfrac{\sigma_1 \eta^i \eta^j}{1 - \sigma_1 \gamma}\Bigr) - \dfrac{\sigma_2 \Bigl(b^i + \dfrac{\sigma_1 \epsilon \eta^i}{(1 - \sigma_1 \gamma)}\Bigr)\Bigl(b^j + \dfrac{\sigma_1 \epsilon \eta^j}{(1 - \sigma_1 \gamma)}\Bigr)}{\tau} \\
                                        &=  a^{ji} + \Bigl(\dfrac{\sigma_1}{1 - \sigma_1 \gamma} - \dfrac{\sigma_1^2 \sigma_2 \epsilon^2}{\tau (1 - \sigma_1 \gamma)^2}\Bigr)\eta^i\eta^j +    \dfrac{\sigma_1 \sigma_2 \epsilon}{\tau (1 - \sigma_1 \gamma)}(b^i\eta^j + b^j\eta^i) - \dfrac{\sigma_2}{\tau}b^ib^j,
                                        \end{split}
                                       \end{equation}
                                       where
                                       \begin{equation*}
                                       \tau = 1 + \sigma_2 \Bigl( w + \dfrac{\sigma_1 \epsilon^2}{(1 - \sigma_1 \gamma)}\Bigr),
                                       \end{equation*}
                                        and
                                        \begin{equation}
                                        det(H_{ij}) = \Bigl[1 + \sigma_2 \Bigl(\omega + \dfrac{\sigma_1 \epsilon^2}{1 - \sigma_1 \gamma}\Bigr)\Bigr](1 - \sigma_1 \gamma)det(a_{ij}). 
                                        \end{equation}
                                        \textbf{Step3.} 
                           
                           Now, we consider
                           \begin{equation}
                           Q_{ji} = a_{ij} - p_1 l_i l_j + p_2 b_i b_j,
                           \end{equation}
                           and
                           \begin{equation}
                            c_i = \sqrt{p_3}\eta_i.
                           \end{equation}
                           By applying Proposition 3.3, we obtain this time:
                           \begin{eqnarray}
                           Q^{ji} &=& a^{ji} + \Bigl(\dfrac{\sigma_1}{(1 - \sigma_1 \gamma)} - \dfrac{\sigma_1^2 \sigma_2 \epsilon^2}{\tau(1 - \sigma_1 \gamma)^2}\Bigr)\eta^i \eta^j \nonumber \\
                           &&
                           - \dfrac{\sigma_1 \sigma_2 \epsilon}{\tau (1 - \sigma_1 \gamma)}(b^i \eta^j + b^j \eta^i) - \dfrac{\sigma_2 b^i b^j}{\tau}.
                           \end{eqnarray}
                           
                                        and
                                         $C_i = P \eta^i + Q b^i$.   
                                         \begin{eqnarray}
                                         P =  \Bigl[ 1 + \Bigl(\dfrac{\sigma_1}{1 - \sigma_1 \gamma} - \dfrac{\sigma_1^2 \sigma_2 \epsilon^2}{\tau (1 - \sigma_1 \gamma)^2}\Bigr)\Bigr]\gamma - \dfrac{\sigma_1 \sigma_2 \epsilon}{\tau(1 - \sigma_1 \gamma)^3}.
                                                                                  \end{eqnarray}
                                          \begin{eqnarray}
                                          Q = -\dfrac{\sigma_2 \epsilon}{\tau} - \dfrac{\sigma_1 \sigma_2 \epsilon \gamma}{\tau (1 - \sigma_1 \gamma)},
                                          \end{eqnarray}
                                          and\\
                                          $C^2 = (P \gamma + Q \epsilon)\sqrt{\sigma_3}$,\\
                                          $ 1 + C^2 = 1 + (P \gamma + Q \epsilon) \sqrt{\sigma_3} \neq 0,$ \\
                                          where
                                        \\
                                        \begin{equation*}
                                        C^i C^j = P^2 \eta^i \eta^j + PQ (\eta^i b^j + \eta^j b^i) + Q^2 b^i b^j.
                                        \end{equation*}
                                        So, the matrix
                                        \begin{equation*}
                                        H_{ij} = a_{ij} - \sigma_1 l_i l_j + \sigma_2 b_i b_j + \sigma_3 \eta_i \eta_j
                                        \end{equation*}
                                        is invertible with
                                        \begin{eqnarray}
                                        H^{ij} &=& a^{ji} + \Bigl[ \dfrac{\sigma_1}{1 - \gamma \sigma_1} - \dfrac{\epsilon^2 \sigma_1^2 \sigma_2}{\tau (1 - \gamma \sigma_1)^2}\Bigr]\eta^i \eta^j - \dfrac{\sigma_2 b^i b^j}{\tau} 
                                        -\dfrac{\epsilon \sigma_1 \sigma_2 (b^i \eta^j + b^j \eta^i)}{\tau (1 - \gamma \sigma_1)} 
                                        \nonumber \\
                                        &&
                                        - \dfrac{P^2 \eta^i \eta^j + PQ(\eta^i b^j + \eta^j b^i) + Q^2 b^i b^j}{1 + (P \gamma + Q \epsilon)\sqrt{\sigma_3}},
                                        \end{eqnarray}
                                        and
                                        \begin{multline}
                                        det(a_{ij} - \sigma_1 l_i l_j + \sigma_2 b_i b_j + \sigma_3 \eta_i \eta_j) =
                                  \Bigl[ 1 + (P \gamma + Q \epsilon)\sqrt{\sigma_3}\Bigr]\Bigl 
                                  [ 1 + \omega + \dfrac{\sigma_1 \epsilon^2}{1 - \sigma_1 \gamma}\Bigr](1 - \sigma_1 \gamma)det(a_{ij}).
                                        \end{multline}
                                     But $g_{ij} = \rho_0H_{ij}$, with $H_{ij} $ from last step. Thus $g^{ji}= \dfrac{1}{\rho_0}H^{ji}$. Using these we immediately obtain the result.
                                         \begin{proposition}
                                          In a non-Hermitian $\mathbb{R}$ - complex Finsler space with special ($\alpha, \beta$)-metric $\dfrac{\beta^2}{\beta - \alpha}$, we have the following properties:
                                          \begin{equation}
                                          \gamma + \bar{\gamma} = l_i \eta^j + l_{\bar{j}} \eta^{\bar{j}} = a_{ij}\eta^j \eta^i + a_{\bar{j}\bar{k}} \eta^{\bar{k}} \eta^{\bar{j}} = 2 \alpha^2.
                                          \end{equation}
                                          \begin{equation}
                                          \theta + \bar{\theta} = b_j \eta^j + b_{\bar{j}} \eta^{\bar{j}} = 2\beta, \delta = \theta.
                                          \end{equation}
                                         \end{proposition}
                                         \textbf{Example.}\\
                                         Consider M = $\mathbb{C}^3$, We set the  metric 
                                         \begin{equation}
                                         \alpha^2 = \exp{z^1 + \bar z^1}|\eta|^2 + \exp(z^2 + \bar{z}^2)|\eta^2|^2 + \exp(z^1 + z^2+z^3+\bar{z}^3)|\eta^3|^2,
                                         \end{equation}  
                                         \begin{equation}
                                         \beta = \dfrac{1}{2}(\exp^{z^2}\eta^2 + \exp^{\bar{z}^2}\eta^2),
                                         \end{equation}
                                         using these values we can form following Matsumoto metric\\
                                         \begin{equation}
                                         F = \dfrac{\exp{z^1 + \bar z^1}|\eta|^2 + \exp(z^2 + \bar{z}^2)|\eta^2|^2 + \exp(z^1 + z^2+z^3+\bar{z}^3)|\eta^3|^2}{\exp{z^1 + \bar z^1}|\eta|^2 + \exp(z^2 + \bar{z}^2)|\eta^2|^2 + \exp(z^1 + z^2+z^3+\bar{z}^3)|\eta^3|^2 - \dfrac{1}{2}(\exp^{z^2}\eta^2 + \exp^{\bar{z}^2}\eta^2)}.
                                         \end{equation}  
                                        	$det(g_{ij}) = \Biggl(\dfrac{\alpha^2(\alpha - 2\beta)}{(\alpha - \beta)^3}\Biggr)^n det(H_{ij})$.

\end{document}